\RequirePackage{fix-cm}
\documentclass[smallextended]{svjour3}       
\smartqed  
\usepackage{graphicx}
\usepackage{subfigure}
\usepackage{multirow}
\usepackage{graphicx}
\usepackage{algorithm}
\usepackage{algorithmic}
\usepackage[a4paper,hmargin={2.54cm,2.54cm},vmargin={3.17cm,3.17cm}]{geometry}

\usepackage{amsmath,amssymb,amsthm,cases}
\usepackage{cite}
\usepackage{enumerate}

\newtheorem{assumption}{Assumption}

\begin{document}

\title{Inertial Proximal Incremental Aggregated Gradient Method \thanks{We are grateful for the support from the National Science Foundation of China (No.11501569, No.61603322).}
}


\author{Xiaoya Zhang        \and
             Wei Peng      \and
             Hui Zhang \and
             Wei Zhu 
}


\institute{Xiaoya Zhang, Wei Peng, Hui Zhang \at
              Department of Mathematics, National University of Defense Technology, Changsha, 410073, Hunan, China. \\
              Tel.: +86-15200860174\\
              \email{zhangxiaoya09@nudt.edu.cn, weipeng0098@126.com, h.zhang1984@163.com.}           
           \and
              Wei Zhu \at
              Shool of Mathematics and Computational Science, Xiangtan University, Xiangtan, 411105, Hunan, China.\\
              \email{zhuwei@xtu.edu.cn}
}

\date{Received: date / Accepted: date}

\maketitle

\begin{abstract}
In this paper, we introduce an inertial version of the Proximal Incremental Aggregated Gradient method(PIAG) for minimizing the sum of smooth convex component functions and a possibly nonsmooth convex regularization function. 
Theoretically, we show that the inertial Proximal Incremental Aggregated Gradient(iPIAG) method enjoys a global linear convergence under a quadratic growth condition, which is strictly weaker than strong convexity, provided that the step size is not larger than some constant. Moreover, we present two numerical experiments which demonstrate that iPIAG outperforms the original PIAG.
\keywords{linear convergence \and inertial method \and quadratic growth condition \and  incremental aggregated gradient \and Lyapunov function\and proximal operator}
\subclass{90C30 \and 90C26 \and 47N10}
\end{abstract}

\section{Introduction}
In this paper, we consider the following composite optimization problem
\begin{eqnarray}\label{M:Pro1}
\text{minimize}_{~x \in \mathbb{R}^d~} \Phi(x) = F\left(x\right) + h\left(x\right)
\end{eqnarray}
where $F\left(x\right):= \sum_{i=1}^N f_i\left(x\right)$ over $N$ training samples describe the fitness to data and 
$h(x)$ is possibly a nonsmooth convex regularization function.  
This problem has a wide range of applications, such as the logistic regression\cite{peng2002introduction, meier2008group} of classification in data mining, the group lasso problem\cite{simon2013sparse} in gene expression, the $\ell_1$-regularized linear least squares problems arising in compressed sensing\cite{chretien2010alternating, hale2007fixed}.

\subsection{The PIAG Method}
\begin{figure*}[htbp!]
 \centering
 \subfigure{
  \includegraphics[width=0.8\textwidth]{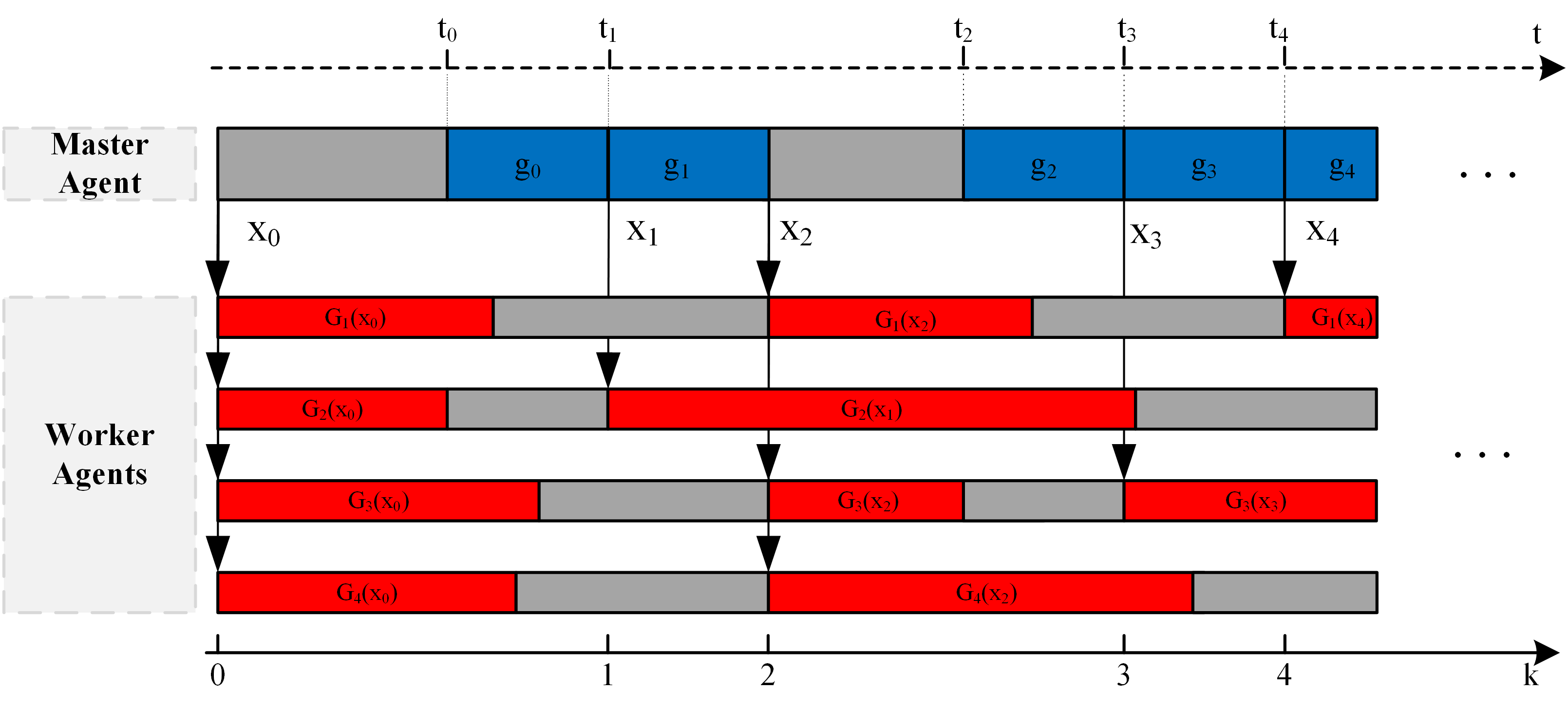}}
 \caption{Process of PIAG with information exchange between master agent and worker agents.\label{Fig:1} In this given example, we consider four worker agents($W=4$) and component functions' gradients are divided into four blocks, that is $\bigcup_{w=1}^4 \mathcal{N}_w =\mathcal{N}$ and $\sum_{w=1}^4|\mathcal{N}_w|=N$. The activation periods of master agent are depicted in blue color while those of worker agents in red. Grey blocks correspond to inactivity of these nodes. Points of the top axis $t$ indicate reading gradients of worker agents by master; Points of the bottom axis $t$ are correponding to returning variable from master agent to worker agents. With the example given in the figure, $g_4 = G_1(x_2) + G_2(x_1) + G_3(x_2) + G_4(x_0)= \sum_{i \in \mathcal{N}_1}\nabla f_i(x_2)+ \sum_{i \in \mathcal{N}_2}\nabla f_i(x_1)+ \sum_{i \in \mathcal{N}_3}\nabla f_i(x_2)+ \sum_{i \in \mathcal{N}_4}\nabla f_i(x_0)$, but not $\sum_{w=1}^4G_w(x_4)$. 
 }
\end{figure*}

Since our work is an inertial version of the Proximal Incremental Aggregated Gradient(PIAG) method, we first introduce the original PIAG, along with a description of asynchronous parallel mechanism.
The asynchronous PIAG algorithm\cite{aytekin2016analysis} for this problem (\ref{M:Pro1}) could be described as:
\begin{align*}
g_k &=\sum_{n=1}^N \nabla f_n(x_{k-\tau_{k}^n}), \tag{PIAG-1}\\
x_{k+1} &= \arg \min_{x}\left\{h(x)+\frac{1}{2\alpha}\|x-(x_k-\alpha g_k)\|^2\right\}.\tag{PIAG-2}
\end{align*}
where $\tau_k^n$ are nonnegative integers. PIAG method makes use of the sum of the most recently evaluated gradients of all component functions  $g_k$ rather than $\sum_{n=1}^N\nabla f_n(x_k)$. 
It can be explained from the perspective of asynchronous parallel mechanism. Actually, in the first step, each worker agent computes gradient information of the component functions on this block and then sends the information to the master agent. In the second step, once the master agent receives the information from one or more worker agents, it will aggregate the gradient information and compute the proximal step (PIAG-2). The agents will not be idle and the optimization speed will not be determined by the slowest worker. In fact, the gradient information, which is aggregated by master agent as (PIAG-1), is inevitably delayed and old for updating.
Figure \ref{Fig:1} illustrates one specific example to explain how to make use of the sum of the most recently evaluated gradients .

Linear convergence of PIAG has been studied under the assumption of strong convexity in a group of recent papers \cite{feyzmahdavian2014delayed, aytekin2016analysis, vanli2016stronger}.
Very recently, in \cite{zhang2017linear} the author proved a linear convergence result under weaker assumptions, quadratic growth condition, than strong convexity. Later on, Zhang et al.\cite{zhang2017Proximal} extended the version to a unified algorithmic framework, called proximal-like incremental aggregated gradient (PLIAG) method, under Bregman distance growth conditions. 

\subsection{Our Work and Contributions}
Accelerated gradient methods have been the hot spot of convex optimization research. In paper \cite{ochs2014ipiano}, the authors developed an inertial proximal gradient algorithm(iPiano) for problem (\ref{M:Pro1}), where the heavy ball method \cite{polyak1964some, stathopoulos2017inertial}  was employed for the PIAG algorithm who set $\tau_k^n$ as $0$. In paper \cite{beck2009fast}, Bech and Teboulle extended Nesterov's accelerated gradient method \cite{nesterov2013introductory} to proximal gradient algorithm. Moreover, in paper \cite{gurbuzbalaban2017convergence}, heavy ball method was employed for Incremental Aggregated Gradiend(IAG) algorithm.
Intuitively, we could extend the inertial methods(include heavy ball method and Nesterov-like accelerated method ) to PIAG algorithm.  

In this paper, we extend PIAG to an unified inertial version with a fully description of asynchronous parallel mechanism, named inertial Proximal Incremental Aggregated Gradient(iPIAG).
In this framework, PIAG with heavy ball method and PIAG algorithm with Nesterov-like acceleration are the specific cases by choosing certain parameters. Moreover, the IAG with heavy method in \cite{gurbuzbalaban2017convergence} is also one special case where $h(x)= 0$. 
To the best of our knowledge, this paper is the first work to build a unified inertial accelerated framework for PIAG algorithm. 

Moreover, we give detailed analysis of this inertial version. By adopting an assumption based on quadratic growth condition, which is strictly weaker than strong convexity, we systematically demonstrate some global linear convergence results for the iPIAG algorithm. In fact, the quadratic growth condition has been proved to be equivalent to other kinds of error bounds under some assumptions in \cite{zhang2016new}. It is possible that theoretical analysis about iPIAG algorithm under other error bounds can be built. But we only focus on study on 
quadratic growth condition and leave others as future work.

Our contributions can be summarized as follow:
\begin{itemize}
\item We introduce the inertial PIAG method, which covers PIAG with momentum and PIAG with Nesterov-like acceleration etc. Besides, we describe the asynchronous parallel parameter server in detail. 
\item We show global linear convergence of iPIAG under weaker assumptions that were used for PIAG in \cite{zhang2017linear,zhang2017Proximal}, by constructing a Lyapunov function. 
\item Experimental results demonstrate that the proposed algorithm is remarkablely faster than the existing PIAG method. 
\end{itemize} 

The rest of the paper is organized as follows. The rest of Section \ref{Sec:1sub2} introduces some basic notations. In Section \ref{Sec:2} we present the iPIAG algorithm for solving problem (\ref{M:Pro1}) in an asynchronous parallel system. In Section \ref{Sec:3sub1}, we first state the assumptions for problem (\ref{M:Pro1}); Then we introduce a lemma in 
\cite{aytekin2016analysis} and on the basis we propose a new modified lemma which could be used in the proof of the main theorem, besides a descent lemma is given; At last, we state the main result of this paper in Section \ref{Sec:3sub2}. Through numerical experiments, Section \ref{Sec:4} verifies that inertial acceleration is conducive to improvement of computing efficiency. Section \ref{Sec:5} concludes the paper.

\subsection*{Notations}\label{Sec:1sub2}
In this paper, we make the following notations:
\begin{itemize}
\item The set of natural numbers $\mathbb{N}=\{1,2,3,\cdots\}$ and the set of natural numbers including zero $\mathbb{N}_0=\{0,1,2,3,\cdots\}$;
\item $\mathcal{N}=\{1,2,\dots,N\}$ and the cardinality of a set $\mathcal{N}$ is $|\mathcal{N}|$;
\item Let $\mathbb{R}^d$ be the Euclidean space with inner product $\langle \cdot, \cdot \rangle$ and induced norm $\|\cdot\|$;
\item For any $x \in \mathbb{R}^d$ and any nonempty set $\Omega \subset \mathbb{R}^d$, the Euclidean distance from $x$ to $\Omega$ is defined by
$d(x,\Omega)= \inf_{y \in \Omega} \|x-y\|$;
\item We let $\mathcal{X}$ be the optimal solution set of problem (\ref{M:Pro1}), and $\Phi^{\ast}$be the associated optimal function value. We always assume that $\mathcal{X}$ is nonempty.
\end{itemize}
\begin{definition}[Proximal Operator\cite{parikh2014proximal}]
Let $f:\mathbb{R}^d \rightarrow \mathbb{R} \cup \{+\infty\}$ be a closed proper convex function. The proximal operator of the scaled function $\alpha f$, where $\alpha >0$, is the function $\mathbf{prox}_{\alpha f}: \mathbb{R}^d \rightarrow \mathbb{R}^d$ defined by
\begin{align}
\mathbf{prox}_{\alpha f} (v)= \arg\min_{x} \left\{f(x) + \frac{1}{2\alpha}\|x-v\|^2\right\}.
\end{align}
\end{definition}

\section{Inertial Proximal Incremental Aggregated Gradient Methods}\label{Sec:2}
Now we introduce the iPIAG algorithm framework. Update of the decision variables $\{x_k,y_k,z_k\}_{k \in \mathbb{N}_0}$ can be computed with setting a constant step size $\alpha>0$, inertial parameters $\eta_1,\eta_2$ and delay parameters $\{\tau^n_k\}_{k \in \mathbb{N}_0}^{n \in \mathcal{N}}$, as 
\begin{equation*}
 \left.\begin{aligned}
        & y_{k+1}=x_k + \eta_1 \left( x_k-x_{k-1} \right), \eta_1 \in [0,1], \\
        & z_{k+1}=   \arg \min_{z}\{h(z)+\frac{1}{2\alpha}\|z-(y_{k+1}-\alpha g_k)\|^2\},\\
        & x_{k+1}=z_{k+1} + \eta_2 \left( z_{k+1}-z_{k} \right), \eta_2 \in [0,1]. 
        \end{aligned}
 \right\}
\qquad 
\eqno{ \text{(iPIAG)}}
\end{equation*}
where $g_k$ is the same with (PIAG-1).
The second step can be represented as a proximal operator step $z_{k+1}=\mathbf{prox}_{\alpha h}\left(y_{k+1} - \alpha g_k\right)$. 
The whole detailed algorithm is described in Algorithm \ref{Alg:1} and Algorithm \ref{Alg:2}, which are corresponding to master agent procedure and worker agent procedures respectively.  
 
\begin{algorithm}[htb]
\caption{Master agent procedure of iPIAG \label{Alg:1}}
\begin{algorithmic}
\STATE{\textbf{Data:} $G_w$ from each worker $w \in \mathcal{W}:=\{1,2,\dots, W\}$}
\STATE{\textbf{Input:} $\alpha$, $\eta_1,\eta_2$, $x_{-1}=x_0=z_0$, regularization function $h(\cdot)$, $K>0$}

\STATE{\textbf{Initialization:} $k=0$}
\STATE{ }

\STATE{1: \textbf{while} $k < K$ \textbf{do}}
\STATE{2: \qquad wait until a set $\mathcal{R}$ of worker agents return their gradients} 
\STATE{3: \qquad \textbf{for all} $w \in \mathcal{W}$ \textbf{do}}
\STATE{4: \qquad \qquad \textbf{if} $w \in \mathcal{R}$}
\STATE{5: \qquad \qquad \qquad Update $G_{w} = \sum_{i \in \mathcal{N}_w} \nabla f_i(x_{k-\tau_k^w}) $}
\STATE{6: \qquad \qquad \textbf{else}}
\STATE{7: \qquad \qquad \qquad Keep old $G_{w}$}
\STATE{8: \qquad \qquad \textbf{end if}}
\STATE{9: \qquad \textbf{end for}}

\STATE{10:\qquad Aggregate the incremental gradients $g_k = \sum_{w \in \mathcal{W}} G_w$}, compute: $x_{k+1}, y_{k+1}, z_{k+1}$ as (iPIAG)

\STATE{11:\qquad \textbf{for all} $w \in \mathcal{R}$ \textbf{do}}
\STATE{12:\qquad \qquad Send $x_{k+1}$ to worker agents $w$}
\STATE{13:\qquad \textbf{end for}}

\STATE{14:\qquad increment $k \leftarrow k+1$}
\STATE{15: \textbf{end while}}
\STATE{ }
\STATE{\textbf{Output} $x_K, y_K, z_K$}
\end{algorithmic}
\end{algorithm}

\begin{algorithm}[htb]
\caption{Procedure for each worker agent $w$ in iPIAG \label{Alg:2}}
\begin{algorithmic}
\STATE{\textbf{Data:}} $x$ and loss functions $\{f_i(\cdot) : i \in \mathcal{N}_w\}$, where $\bigcup_{w \in \mathcal{W}} \mathcal{N}_w = \mathcal{N}$ and $\mathcal{N}_{w_1} \bigcap \mathcal{N}_{w_2} =\emptyset(\forall w_1,w_2\in \mathcal{W}, w_1 \neq w_2 $)

\STATE{1: \textbf{repeat}}
\STATE{2: \qquad Receive $x \leftarrow x_{k}$ from master}
\STATE{3: \qquad Calculate: $G_w = \sum_{i \in \mathcal{N}_w} \nabla f_i(x)$ }
\STATE{4: \qquad Send $G_w$ to master agent with a delay of $\tau_k^w$}
\STATE{5: \textbf{until} EXIT received}
\end{algorithmic}
\end{algorithm}

Since iPIAG is an unified inertial extension of PIAG, it is worth mentioning that several recent works are the special cases.
\begin{enumerate}[(i)]
    \item When $\eta_1=\eta_2=0$, iPIAG reduces to PIAG;
	\item When $h(x) \equiv 0$, and $\eta_2=0$, iPIAG reduces to the incremental aggregated gradient algorithm with momentum(IAG with momentum) in paper \cite{gurbuzbalaban2017convergence};
	\item When $\{\tau^n_k\}_{k \in \mathbb{N}_0}^{n \in \mathcal{N}}=\{0\}$, and $\eta_2=0$, iPIAG reduces to the ipiasco algorithm\cite{ochs2015ipiasco};
	\item When $\eta_2 = 0$, iPIAG reduces to a PIAG with heavy ball method, abbreviated as PIAG-M;
	\item When $\eta_1 = 0$, iPIAG reduces to a PIAG with Nesterov-like acceleration, abbreviated as PIAG-NeL. 
\end{enumerate}
\begin{remark}
The above (iii) and (iv) is also first proposed in our paper. We could give concrete parameter selection rules in Corollary \ref{Corol:1} and \ref{Corol:2} for PIAG-M and PIAG-NeL respectively to guarantee linear convergence rate.
\end{remark}

\section{Linear Convergence of iPIAG Method}\label{Sec:3}
In this section, we show that iPIAG converges linearly. We first make some assumptions and then prove the main results through Theorem \ref{Th:1}.

\subsection{Assumptions}\label{Sec:3sub1}
We need the following assumptions.
\begin{assumption} \label{Assum:1}
Each component function $f_n$ is convex with $L_n$-continuous gradient:
$$\|\nabla f_n(x) - \nabla f_n(y) \| \leqslant L_n\|x-y\|, ~~~~\forall x,y \in \mathbb{R}^d.$$
Moreover, we can get the conclusion that the sum function $F$ is convex with $L$-continuous gradient where $L=\sum_{n=1}^N L_n$.
\end{assumption} 

\begin{assumption} \label{Assum:2}
The regularization function $h: \mathbb{R}^d \rightarrow (-\infty, \infty ]$ is proper, closed, convex and subdifferentiable everywhere in its effective domain, that is $\partial h(x) \neq \emptyset$ for all $x \in \{y \in \mathbb{R}^d: h(y)< \infty\}$. 
\end{assumption} 

\begin{assumption}\label{Assum:3}
The time-varying delays $\tau_k^n$ are bounded, i.e., there is some integer $\tau \in \mathbb{N}_0$ such that
$$\tau_k^n \in \{0,1,\dots ,\tau\}, ~~~~\forall k \in \mathbb{N}_0, n \in \mathcal{N}$$
holds. Such $\tau$ is called the delay parameter.
\end{assumption} 
\begin{assumption}\label{Assum:4}
The objective function  $\Phi$ satisfies the quadratic growth condition, meaning there is a real number  $\beta > 0$ such that
$$ \Phi(x)-\Phi^* \geqslant \frac{\beta}{2}d^2(x, \mathcal{X}), ~~~~x \in \mathbb{R}^d.$$
\end{assumption} 
Assumptions \ref{Assum:1}-\ref{Assum:3} are standard prerequisites for proving linear convergence of PIAG and IAG in papers \cite{gurbuzbalaban2017convergence, aytekin2016analysis, vanli2016stronger, vanli2016global}.
The difference between our assumptions and what were made in the above papers lies in Assumption \ref{Assum:4}: strong convexity of $\Phi$ is replaced. Assumption \ref{Assum:4} was recently used in \cite{zhang2017linear, zhang2017Proximal} for PIAG and it has been shown as a weaker condition than strong convexity for problem (\ref{M:Pro1}).

\subsection{Key Lemmas and Main Results}\label{Sec:3sub2}
\subsubsection{Key Lemmas}
The following lemma first appeared in \cite{aytekin2016analysis} and was employed as an important tool for linear convergence analysis of PIAG in \cite{vanli2016stronger, zhang2017linear, zhang2017Proximal} recently.
\begin{lemma}\label{Lem:1}
Assume that the nonnegative sequences $\{V_k\}$ and $\{w_k\}$ satisfy the following inequality:
\begin{eqnarray}
V_{k+1} \leqslant aV_k -b w_k +c\sum_{j=k-k_0}^{k}w_j, ~~~~k \geq 0
\end{eqnarray}
for some real numbers a $a \in (0, 1)$ and $b, c \geqslant 0$, and some nonnegative integer $k_0$. Assume that $w_k = 0$ for $k < 0$, and the following holds:
\begin{eqnarray}
\frac{c}{1-a} \frac{1-a^{k_0+1}}{a^{k_0}} \leqslant b.
\end{eqnarray}
Then, $V_k \leqslant a^kV_0$ for all $k\geqslant 0$.
\end{lemma}

In order to show the linear convergence of iPIAG, we need a slightly complicated result, which bears the spirit of Lemma \ref{Lem:1}.
\begin{lemma}\label{Lem:2}
Assume that the non-negative sequences $\{V_k\}$ and $\{\omega_k\}$ satisfy the following inequality for $k\geq 1$:
\begin{align}\label{ineq1}
V_{k+1}\leq AV_k+BV_{k-1}-b_1\omega_k+b_2\omega_{k-1}+c\sum_{j=k-k_0}^k\omega_j, ~~~~k \geq 0
\end{align}
for positive real numbers $A,B$ satisfying $A+B<1$ and $b_1,b_2,c\geq 0$, and some integer $k_0\in\mathbb{N}$. Assume that $\omega_k=0$ for $k<0$.
Take constants $\alpha_1,\alpha_2,a$ to satisfy
\begin{equation}\label{cond5}
\left \{
\begin{aligned} 
\alpha_1+\alpha_2=1\\
A=\alpha_1a\\
B=\alpha_2a^2\\
0<a<1\\
\alpha_1,\alpha_2\geq 0.
\end{aligned}
\right.
\end{equation}
If the following holds:
\begin{align}\label{cond1}
\frac{c}{1-a}\frac{1-a^{k_0+1}}{a^{k_0}}\leq b_1-\frac{b_2}{a},
\end{align}
then $V_k$ is linear convergent in the sense that
\begin{align*}
V_{k}\leq a^{k-1}\left(V_1+ aV_0+b_1\omega_0\right),~~~~\forall k\geq 1.
\end{align*}
\begin{proof}
Note that $A=\alpha_1 a$ and $B=\alpha_2 a^2$ in (\ref{cond5}),
the inequality (\ref{ineq1}) can be rewritten as:
\begin{align}\label{Lem2:proof1}
V_{k+1}\leq \alpha_1aV_k+\alpha_2a^2V_{k-1}-b_1\omega_k+b_2\omega_{k-1}+c\sum_{j=k-k_0}^k\omega_j.
\end{align}

By dividing both sides of (\ref{Lem2:proof1}) by $a^{k+1}$ and summing the resulting inequality up from $k=1$ to $K(K \geq 1)$, we derive that
\begin{flalign}
\sum_{k=1}^K \frac{V_{k+1}}{a^{k+1}} \leq &\alpha_1\sum_{k=1}^K \frac{V_{k}}{a^k}+\alpha_2\sum_{k=0}^{K-1}\frac{V_{k}}{a^{k}}-b_1\sum_{k=1}^K\frac{\omega_k}{a^{k+1}}
+\frac{b_2}{a}\sum_{k=0}^{K-1}\frac{\omega_{k}}{a^{k+1}}+\sum_{k=1}^{K}\left[\frac{c}{a^{k+1}}\sum_{j=k-k_0}^k\omega_j\right] \nonumber\\
\leq &\sum_{k=0}^K\frac{V_k}{a^k}-\left(b_1-\frac{b_2}{a}\right)\sum_{k=0}^K\frac{\omega_k}{a^{k+1}} 
+\sum_{k=0}^{K}\left[\frac{c}{a^{k+1}}\sum_{j=k-k_0}^k\omega_j\right]+\frac{b_1\omega_0}{a} \label{Eq:Lem2_2}
\end{flalign}

Since $w_k=0(k<0), w_k \geq 0(k \geq 0)$ and $a>0$, thus we get
\begin{flalign}
\sum_{k=0}^{K}\left[\frac{c}{a^{k+1}}\sum_{j=k-k_0}^k\omega_j\right]
&\leq \frac{c}{a}\left(\frac{1}{a^{-k_0}} + \frac{1}{a^{-k_0+1}}+ \cdots + \frac{1}{a^{0}} \right)w_{-k_0}+ \frac{c}{a}\left(\frac{1}{a^{-k_0+1}} + \frac{1}{a^{-k_0+2}} + \cdots + \frac{1}{a^{1}} \right)w_{-k_0+1} \nonumber \\
&~~~~ + \cdots +  \frac{c}{a}\left(\frac{1}{a^{-k_0+K}} + \frac{1}{a^{-k_0+K+1}}+ \cdots + \frac{1}{a^{K}} \right)w_{-k_0+K}\nonumber\\
&~~~~ + \cdots +  \frac{c}{a}\left(\frac{1}{a^{K}} + \frac{1}{a^{K+1}} + \cdots + \frac{1}{a^{K+k_0}} \right)w_{K}\nonumber\\
&\leq \sum_{j=0}^{K+k_0}\frac{c}{a}\left(\frac{1}{a^{j-k_0}}+\dots + \frac{1}{a^{j}}\right)w_{j-k_0}\nonumber\\
&= \sum_{j=0}^{K+k_0}\frac{c}{a}\left(\frac{1}{a^{j-k_0}}+\dots + \frac{1}{a^{j}}\right)w_{j-k_0}\nonumber\\
&= \sum_{k=-k_0}^{K}\frac{c}{a}\left(\frac{1}{a^{k}}+\dots + \frac{1}{a^{k+k_0}}\right)w_{k}\nonumber\\
&= \sum_{k=0}^{K}\frac{c}{a}\left(\frac{1}{a^{k}}+\dots + \frac{1}{a^{k+k_0}}\right)w_{k}\label{Eq:Lem2_3}.
\end{flalign}
Combining (\ref{Eq:Lem2_3}) with (\ref{Eq:Lem2_2}), yields
\begin{flalign*}\sum_{k=1}^K \frac{V_{k+1}}{a^{k+1}} &\leq
\sum_{k=0}^K\frac{V_k}{a^k}+ \frac{b_1\omega_0}{a}+\sum_{k=0}^K\left[c\left(1+\frac{1}{a}+\cdots+\frac{1}{a^{k_0}}\right)-\left(b_1-\frac{b_2}{a}\right)\right]\frac{\omega_k}{a^{k+1}}\\
&= \sum_{k=0}^K\frac{V_k}{a^k} + \frac{b_1\omega_0}{a}+
\sum_{k=0}^K\left[   \frac{c}{1-a}\frac{1-a^{k_0+1}}{a^{k_0}}-\left(b_1-\frac{b_2}{a}\right)\right]\frac{\omega_k}{a^{k+1}}.
\end{flalign*}
By the condition (\ref{cond1}), we obtain 
\begin{flalign}
\sum_{k=1}^K \frac{V_{k+1}}{a^{k+1}} \leq
\sum_{k=0}^K\frac{V_k}{a^k}+ \frac{b_1\omega_0}{a},
\end{flalign}
that is $\frac{V_{K+1}}{a^{K+1}} \leq V_0 + \frac{V_1}{a}+ \frac{b_1\omega_0}{a}$ for $K \geq 1$. Besides, we know $V_1 \leq V_1+aV_0+b_1\omega_0$.
Therefore, for $\forall K\geq 1$, we have
\begin{align*}
V_{K}\leq a^{K-1}\left(V_1+aV_0+b_1\omega_0\right).
\end{align*}
\end{proof}
\end{lemma}

The following is a descent-type lemma. 
\begin{lemma}\label{Lem:3}
Suppose that the standard Assumptions \ref{Assum:1}-\ref{Assum:3} hold. Let 
$$ \Delta _k^1 :=\frac{L(\tau+2)}{2}\sum_{j=k-\tau-1}^k\|z_{j+1}-z_j\|^2$$
and $$\Delta_k^2 := \frac{\eta_1(1+\eta_2)^2}{\alpha}\sum_{j=k-2}^{k-1}\|z_{j+1}-z_{j}\|^2.$$
Then, we have that
\begin{flalign}\label{Eqn:Lem3}
\Phi(z_{k+1})&\leq \Phi(x)+\frac{1+\eta_2}{2\alpha}\|x-z_k\|^2-\frac{1-\eta_1}{2\alpha}\|x-z_{k+1}\|^2 -\frac{1}{4\alpha}\|z_{k+1}-z_{k}\|^2 \nonumber \\
&~~~~+ \frac{\eta_2 +2\eta_2^2 }{2\alpha}\|z_k-z_{k-1}\|^2 +\Delta_k^1+\Delta_k^2, k\geq 0.
\end{flalign}
\begin{proof}
Since each component function $f_n(x)$ is convex with $L_n$-continuous gradient, we derive that
\begin{flalign}\label{Eqn:Lem3_1}
f_n(z_{k+1}) & \leqslant f_n(x_{k-\tau_k^n}) + \langle \bigtriangledown f_n(x_{k-\tau_k^n}), z_{k+1}- x_{k-\tau_k^n} \rangle + \frac{L_n}{2} \|z_{k+1}- x_{k-\tau_k^n}\|^2 \nonumber \\ 
&\leqslant f_n(x) + \langle \bigtriangledown f_n(x_{k-\tau_k^n}), z_{k+1}- x \rangle + \frac{L_n}{2} \|z_{k+1}- x_{k-\tau_k^n}\|^2, n=1,\cdots, N,
\end{flalign}
where the second inequality follows from the convexity of $f_n(x)$. Summing (\ref{Eqn:Lem3_1})
over all components functions and using the expression of $g_k$, yield
\begin{eqnarray}\label{Eqn:Lem3_2}
F(z_{k+1}) \leqslant F(x) + \langle g_k, z_{k+1}- x \rangle + \sum_{n=1}^N \frac{L_n}{2} \|z_{k+1}- x_{k-\tau_k^n}\|^2. 
\end{eqnarray}

Noting that $z_{k+1}$ is the minimizer of the $\frac{1}{\alpha}$-strongly convex function:
\begin{eqnarray*}
z \rightarrow h(z) + \frac{1}{2\alpha}\| z-(y_{k+1} -\alpha g_k)\|^2, 
\end{eqnarray*}
for all $x \in \mathbb{R}^d$, we have 
\begin{flalign}\label{Eqn:Lem3_3}
h(z_{k+1}) + \frac{1}{2\alpha}\| z_{k+1}-(y_{k+1} -\alpha g_k)\|^2 \leq h(x)  
+ \frac{1}{2\alpha}\| x-(y_{k+1} -\alpha g_k)\|^2-\frac{1}{2\alpha}\|x-z_{k+1}\|^2.
\end{flalign}
After rearranging the terms of (\ref{Eqn:Lem3_3}), we further get
\begin{flalign}\label{Eqn:Lem3_4}
\langle z_{k+1}-x, g_k \rangle  \leq &h(x) - h(z_{k+1})+ \frac{1}{2\alpha}\| x- y_{k+1}\|^2-\frac{1}{2\alpha}\| z_{k+1}-x\|^2 - \frac{1}{2\alpha}\| z_{k+1}-y_{k+1}\|^2, \forall x \in \mathbb{R}^d.
\end{flalign}
Combining (\ref{Eqn:Lem3_2}) and (\ref{Eqn:Lem3_4}), and noting that $y_{k+1}=x_k +\eta_1(x_k-x_{k-1})$, we derive that 
\begin{flalign}\label{Eqn:Lem3_5}
\Phi(z_{k+1}) &\leq \Phi(x) + \frac{1}{2\alpha}\| x- y_{k+1}\|^2  - \frac{1}{2\alpha}\| z_{k+1}-x\|^2 -\frac{1}{2\alpha}\| z_{k+1}-y_{k+1}\|^2 + \sum_{n=1}^N \frac{L_n}{2} \|z_{k+1}- x_{k-\tau_k^n}\|^2 \nonumber\\
&= \Phi(x) + \frac{1}{2\alpha}\| x- x_k\|^2  -\frac{1}{2\alpha}\| z_{k+1}-x_{k}\|^2 - \frac{1}{\alpha} \langle x-x_k, \eta_1(x_k-x_{k-1})\rangle + \frac{1}{\alpha} \langle z_{k+1}-x_k, \eta_1(x_k-x_{k-1})\rangle\nonumber\\
&~~~~ - \frac{1}{2\alpha}\| z_{k+1}-x\|^2  + \sum_{n=1}^N \frac{L_n}{2} \|z_{k+1}- x_{k-\tau_k^n}\|^2 \nonumber\\
&\leq \Phi(x) + \frac{1}{2\alpha}\| x- x_{k}\|^2 - \frac{1}{2\alpha}\| z_{k+1}-x_k\|^2 - \frac{1}{2\alpha}\| z_{k+1}-x\|^2 + \frac{\eta_1}{2\alpha}(\|x-z_{k+1}\|^2+\|x_{k}-x_{k-1}\|^2)\nonumber \\
&~~~~+ \sum_{n=1}^N \frac{L_n}{2} \|z_{k+1}- x_{k-\tau_k^n}\|^2.
\end{flalign}

Now we bound the three terms $\frac{1}{2\alpha}\|x-x_k\|^2$, $-\frac{1}{2\alpha}\|z_{k+1}-x_k\|^2$ and $\frac{\eta_1}{2\alpha}\|x_{k}-x_{k-1}\|^2$ in (\ref{Eqn:Lem3_5}). Substituting $x_k=z_k+\eta_2(z_k-z_{k-1})$ in (iPIAG) into these terms, we obtain that
\begin{itemize}
\item[$\bullet$] $\frac{1}{2\alpha}\|x-x_k\|^2$ is bounded by 
\begin{flalign}\label{Eqn:Lem3_6}
\frac{1}{2\alpha}\|x-x_k\|^2 &= \frac{1}{2\alpha}\left[(1+\eta_2)^2\left\|\frac{x-z_k}{1+\eta_2} + \frac{\eta_2(z_{k-1}-z_k)}{1+\eta_2}\right\|^2\right] \nonumber\\
&\leq \frac{1}{2\alpha}\left[(1+\eta_2)\|x-z_k\|^2 + \eta_2(1+\eta_2) \|z_k-z_{k-1}\|^2\right],
\end{flalign}
where the last inequality follows from the convexity of function $s(u):=\|u\|^2$. 
\item[$\bullet$]   $-\frac{1}{2\alpha}\|z_{k+1}-x_k\|^2$ is bounded by
\begin{flalign}\label{Eqn:Lem3_7}
-\frac{1}{2\alpha}\|z_{k+1}-x_k\|^2 &= -\frac{1}{2\alpha}(\|z_{k+1}- \left(z_k+\eta_2(z_k-z_{k-1})\right)\|^2)\nonumber\\
&\leq -\frac{1}{2\alpha}\left(\frac{1}{2}\|z_{k+1}-z_{k}\|^2 - \eta_2^2\|z_{k}-z_{k-1}\|^2\right),
\end{flalign}
where the last inequality follows from the fact that
\begin{flalign*}
\|z_{k+1}-z_{k}\|^2  \leq 2\left( \|z_{k+1}- \left(z_k+\eta_2(z_k-z_{k-1})\right)\|^2 + \eta_2^2\|z_{k}-z_{k-1}\|^2 \right).
\end{flalign*}
\item[$\bullet$] $\frac{\eta_1}{2\alpha}\|x_{k}-x_{k-1}\|^2$ is bounded by
\begin{flalign}\label{Eqn:Lem3_8}
\frac{\eta_1}{2\alpha}\|x_{k}-x_{k-1}\|^2 &= \frac{\eta_1}{2\alpha}(\|z_k+\eta_2(z_k-z_{k-1})-(z_{k-1}+\eta_2(z_{k-1}-z_{k-2}))\|^2)  \nonumber\\
&\leq \frac{\eta_1(1+\eta_2)^2}{\alpha}\|z_k-z_{k-1}\|^2 + \frac{\eta_1\eta_2^2}{\alpha} \|z_{k-1}-z_{k-2}\|^2,
\end{flalign}
where the last inequality follows from the inequality $\|a+b\|^2 \leq 2(\|a\|^2+\|b\|^2)$.
\end{itemize}
Combining (\ref{Eqn:Lem3_6})-(\ref{Eqn:Lem3_8}) with (\ref{Eqn:Lem3_5}), we get
\begin{flalign}\label{Eqn:Lem3_9}
\Phi(z_{k+1}) &\leq
\Phi(x) + \frac{1}{2\alpha}\left[(1+\eta_2)\|x-z_k\|^2 + \eta_2(1+\eta_2) \|z_k-z_{k-1}\|^2\right]
-\frac{1}{2\alpha}\left(\frac{1}{2}\|z_{k+1}-z_{k}\|^2 - \eta_2^2\|z_{k}-z_{k-1}\|^2\right)\nonumber\\
&~~~~- \frac{1}{2\alpha}\| z_{k+1}-x\|^2 + \frac{\eta_1}{2\alpha}\|x-z_{k+1}\|^2+
\frac{\eta_1(1+\eta_2)^2}{\alpha}\|z_k-z_{k-1}\|^2 + \frac{\eta_1\eta_2^2}{\alpha} \|z_{k-1}-z_{k-2}\|^2 \nonumber \\
&~~~~+ \sum_{n=1}^N \frac{L_n}{2} \|z_{k+1}- x_{k-\tau_k^n}\|^2\nonumber \\
&=\Phi(x) + \frac{1+\eta_2}{2\alpha}\|x-z_k\|^2 +  \frac{\eta_2(1+2\eta_2)}{2\alpha}\|z_k-z_{k-1}\|^2 -\frac{1}{4\alpha}\|z_{k+1}-z_{k}\|^2- \frac{1-\eta_1}{2\alpha}\| z_{k+1}-x\|^2 \nonumber \\
&~~~~+\frac{\eta_1(1+\eta_2)^2}{\alpha}\|z_k-z_{k-1}\|^2 + \frac{\eta_1\eta_2^2}{\alpha} \|z_{k-1}-z_{k-2}\|^2 + \sum_{n=1}^N \frac{L_n}{2} \|z_{k+1}- x_{k-\tau_k^n}\|^2
\end{flalign}

Since $x_{k-\tau_k^n} = z_{k-\tau_k^n}+ \eta_2(z_{k-\tau_k^n}- z_{k-\tau_k^n-1})$, we have
\begin{flalign}\label{Eqn:Lem3_10}
\sum_{n=1}^N\frac{L_n}{2}\|z_{k+1}-x_{k-\tau_k^n}\|_2^2
&= \sum_{n=1}^N\frac{L_n}{2}\|z_{k+1}-z_{k}+\cdots+z_{\tau_k^i}-x_{k-\tau_k^n}\|^2 \nonumber \\
&=\sum_{n=1}^N\frac{L_n}{2}\|z_{k+1}-z_{k}+\cdots+\eta_2 z_{k-\tau_k^n}-\eta_2 z_{k-\tau_k^n-1}\|^2\nonumber \\
&\leq\frac{L(\tau+2)}{2}\sum_{j=k-\tau-1}^k\|z_{j+1}-z_{j}\|^2 \equiv \Delta_k^1,
\end{flalign}
where the last inequality follows the Jensen’s inequality. Besides, since $0 \leq \eta_2 \leq 1$, $\eta_2^2 \leq (1+\eta_2)^2$, we have
\begin{flalign}\label{Eqn:Lem3_11}
\frac{\eta_1(1+\eta_2)^2}{\alpha}\|z_k-z_{k-1}\|^2 +\frac{\eta_1\eta_2^2}{\alpha}\|z_{k-1}-z_{k-2}\|^2 \leq \frac{\eta_1(1+\eta_2)^2}{\alpha}\sum_{j=k-2}^{k-1}\|z_{j+1}-z_{j}\|^2 \equiv \Delta_k^2. 
\end{flalign}
By substituting (\ref{Eqn:Lem3_10}) and (\ref{Eqn:Lem3_11}) into (\ref{Eqn:Lem3_9}), the desired result follows.
\end{proof}
\end{lemma}

\subsubsection{Main Results}
Now we state the main result of this paper.  We first define a Lyapunov function
$$\Psi(x) := \Phi(x)-\Phi^* + \frac{1-\eta_1}{2\alpha}d^2(x,\mathcal{X}).$$

\begin{theorem}\label{Th:1}
Suppose that Assumptions \ref{Assum:1}-\ref{Assum:4} hold. Assume that the step-size and inertial parameters satisfy:
\begin{eqnarray*}
\alpha \leqslant \alpha_0 := \frac{ (W+1)^{\frac{1}{\tau+3}}-1}{\beta},
\end{eqnarray*}
where $W = \frac{\beta}{16C_1\beta + 2L(\tau+2)}(C_1 \in [0, \frac{1}{2}))$, and
\begin{flalign*}
&\eta_1 = \min\{C_1\alpha\beta,1\} , \\
&\eta_2 \in \left[0,\min\left\{\frac{\alpha\beta}{2}, \frac{1}{1+\alpha\beta-\eta_1} \left(\frac{1}{4}-\frac{L(\tau+2)\alpha+8\eta_1}{2\alpha\beta}\left((\alpha\beta+1)^{\tau+3}-1\right)
\right)
\right\}\right].
\end{flalign*}
 
Let $\mathbf{C}:= \Psi(z_1)+ \frac{1+\eta_2}{1+\alpha\beta-\eta_1}\Psi(z_0)+\frac{1}{4\alpha}\|z_1-z_0\|^2_2$. Then, $\{\Psi(z_k)\}_{k=1}^{\infty}$ converges linearly in the sense that:
\begin{eqnarray} \label{Equ:ly}
\Psi(z_k)\leq \left(\frac{1+\eta_2}{1+\alpha\beta-\eta_1}\right)^{k}\mathbf{C}, k\geq 1.
\end{eqnarray}
In particular,
\begin{enumerate}[(i)]
\item iPIAG method attains a global linear convergence in function values
\begin{eqnarray}\label{Equ:obj}
\Phi(z_k)-\Phi^* \leq \left(\frac{1+\eta_2}{1+\alpha\beta-\eta_1}\right)^{k}\mathbf{C},k\geq 1;
\end{eqnarray}
\item There exists a global linear convergence in distances of the iterates to the optimal solution set:
\begin{eqnarray}\label{Equ:point}
d^2(z_k,\mathcal{X}) \leqslant \frac{2\alpha}{1-\eta_1}\left(\frac{1+\eta_2}{1+\alpha\beta-\eta_1}\right)^{k}\mathbf{C},k\geq 1.
\end{eqnarray}
\end{enumerate}
\begin{proof}
Since the set $\mathcal{X}$ is convex and nonempty, the projection point of $z$ onto $\mathcal{X}$ is unique. We use $z^{\ast}$ to stand for this projection point. Note that  $\Phi(z_k^{\ast})=\Phi^*$. According to Lemma \ref{Lem:3}, we obtain
\begin{flalign}\label{Eqn:Th_1}
\Phi(z_{k+1}) - \Phi^* + \frac{1-\eta_1}{2\alpha}\|z_k^{\ast}-z_{k+1}\| \leq \frac{1+\eta_2}{2\alpha}\|z_k^{\ast}-z_{k}\|^2 -\frac{1}{4\alpha}\| z_{k+1}-z_{k} \|^2 +\frac{\eta_2 +2\eta_2^2 }{2\alpha}\|z_k-z_{k-1}\|^2 + \Delta_k^2+ \Delta_k^1.
\end{flalign}
Since $z_k^{\ast}\in \mathcal{X}$, following from the definition of projection, it holds that
$$\|z_k^{\ast}-z_{k+1}\|^2 \geq \|z_{k+1}^{\ast}-z_{k+1}\|^2 = d^2(z_{k+1}, \mathcal{X}).$$
Thus by using the expression of the Lyapunov function $\Psi$, we further obtain
\begin{flalign}\label{Eqn:Th_2}
\Psi(z_{k+1}) \leq \frac{1+\eta_2}{2\alpha}\|z_k^{\ast}-z_{k}\|^2 -\frac{1}{4\alpha}\| z_{k+1}-z_{k} \|^2+ \frac{\eta_2 +2\eta_2^2 }{2\alpha}\|z_k-z_{k-1}\|^2+ \Delta_k^2+ \Delta_k^1.
\end{flalign}

By using the quadratic growth condition, we obtain
$$\|z_k^{\ast}-z_k\|^2 = d^2(z_{k}, \mathcal{X}) \leq \frac{2}{\beta}(\Phi(z_{k}) - \Phi^*),$$
and hence
\begin{flalign}\label{Eqn:Th_3}
\|z_k^{\ast}-z_{k}\|^2 \leq p\|z_k^{\ast}-z_{k}\|^2 + \frac{2q}{\beta}(\Phi(z_k)-\Phi^{\ast})
\end{flalign}
with $p+q=1, p, q \geq 0$.
Picking $p=\frac{1-\eta_1}{\alpha\beta+1-\eta_1},q=\frac{\alpha\beta}{\alpha\beta+1-\eta_1}$ and combining (\ref{Eqn:Th_2}) and (\ref{Eqn:Th_3}), we obtain
\begin{align}\label{Eqn:Th_4}
\Psi(z_{k+1}) \leq \frac{1+\eta_2}{\alpha\beta+1-\eta_1}\Psi(z_k) -\frac{1}{4\alpha}\| z_{k+1}-z_{k} \|^2+\frac{\eta_2 +2\eta_2^2 }{2\alpha}\|z_k-z_{k-1}\|^2+ \Delta_k^2+ \Delta_k^1.
\end{align}

We define $a(\eta_1,\eta_2):=\frac{1+\eta_2}{\alpha\beta+1-\eta_1}$. We will need $a(\eta_1,\eta_2) < 1$, which is equivlent to requiring $\eta_1+ \eta_2 < \alpha\beta$.  
To this end, we set 
\begin{flalign}\label{Eqn:eta}
\eta_1=C_1\alpha\beta, \eta_2=C_2\alpha\beta,
\end{flalign}
where $0 \leq C_1 + C_2 < 1, C_1,C_2 \geq 0$.

Now we estimate the convergence rate via discussing two different types of $\tau$.
\begin{itemize}
\item Case $\tau \geq 1$:
Since　$~\eta_2 \leq 1$, we get $ \frac{\eta_1(1+\eta_2)^2}{\alpha} \leq 4C_1\beta$ and further obtain
\begin{align}\label{Eqn:Th_5}
\Delta_k^1 + \Delta_k^2 &= \frac{L(\tau+2)}{2}\sum_{j=k-\tau-1}^k\|z_{j+1}-z_j\|^2
+ \frac{\eta_1(1+\eta_2)^2}{\alpha}\sum_{j=k-2}^{k-1}\|z_{j+1}-z_{j}\|^2\\
&\leq \left(\frac{L(\tau+2)}{2}+ 4C_1\beta\right)\sum_{j=k-\tau-1}^k\|z_{j+1}-z_j\|^2,
\end{align}
In order to apply Lemma \ref{Lem:2}, take
$b_1=\frac{1}{4\alpha}$, $b_2=\frac{\eta_2+2\eta_2^2}{2\alpha}$, $c=\frac{L(\tau+2)}{2}+4C_1\beta$, $k_0=\tau+1$.
From (\ref{Eqn:Th_4}), we obtain
\begin{align}\label{Eqn:Th_6}
\Psi(z_{k+1}) \leq a(\eta_1,\eta_2)\Psi(z_k) - b_1 \|z_{k+1}-z_{k} \|^2+ b_2\|z_{k}-z_{k-1} \|^2+c\sum_{j=k-k_0}^k\|z_{j+1}-z_j\|^2.
\end{align}

To apply Lemma \ref{Lem:2}, we need the parameters to satisfy (\ref{cond1}). First of all, note that
\begin{align}\label{Eqn:Th_7}
\frac{b_2}{a(\eta_1,\eta_2)} 
=\frac{1+(1-C_1)\alpha\beta}{\eta_2+1}\frac{2\eta_2^2+\eta_2}{2\alpha} 
\leq \frac{1+(1-C_1)\alpha\beta}{2\alpha}\frac{2\eta_2^2+2\eta_2}{\eta_2+1}
= \eta_2 \frac{1+(1-C_1)\alpha\beta}{\alpha}.
\end{align}
Second, we denote $r(t)=\frac{c}{1-t}\frac{1-t^{k_0+1}}{t^{k_0}}$. Note that 
\begin{align*}
r(t)=c\cdot\sum_{j=0}^{k_0}\left(\frac{1}{t}\right)^j
\end{align*}
is monotonically decreasing and $a(\eta_1,\eta_2) \geq a(0,0)$. Thus, $r(a(\eta_1,\eta_2)) \leq r(a(0,0))$, which means 
\begin{align}\label{Eqn:Th_8}
b_1-\frac{c}{1-a(0,0)}\frac{1-a(0,0)^{k_0+1}}{a(0,0)^{k_0}}
\leq b_1-\frac{c}{1-a(\eta_1,\eta_2)}\frac{1-a(\eta_1,\eta_2)^{k_0+1}}{a(\eta_1,\eta_2)^{k_0}}.
\end{align}

In light of (\ref{Eqn:Th_7}) and (\ref{Eqn:Th_8}), we have (\ref{cond1}) if the parameters satisfy 
\begin{align}\label{Eqn:Th_9}
\eta_2 \frac{1+(1-C_1)\alpha\beta}{\alpha} 
&\leq b_1-\frac{c}{1-a(0,0)}\frac{1-a(0,0)^{k_0+1}}{a(0,0)^{k_0}}\nonumber\\
& = \frac{1}{4\alpha}-\frac{L(\tau+2)+8C_1\beta}{2\alpha\beta}\left((\alpha\beta+1)^{\tau+2}-1\right).
\end{align}
Thus we could choose $\eta_1,\eta_2$ from the nonempty intervals respectively,
\begin{flalign}
&\mathcal{I}_1= [0, \frac{\alpha\beta}{2}), \text{that is } 0 \leq C_1< \frac{1}{2},\label{Eqn:Th_10}\\
&\mathcal{I}_2=\left[0,\min\left\{\frac{\alpha\beta}{2}, \frac{1}{1+(1-C_1)\alpha\beta} \left(\frac{1}{4}-\frac{L(\tau+2)+8C_1\beta}{2\beta}\left((\alpha\beta+1)^{\tau+2}-1\right)
\right)
\right\}\right].\label{Eqn:Th11}
\end{flalign}
In order to show that $\eta_2$ is selected properly, there must be
\begin{align*}
\frac{1}{4}-\frac{L(\tau+2)+8C_1\beta}{2\beta}\left((\alpha\beta+1)^{\tau+2}-1\right)>0,
\end{align*}
which is equivalent to requiring the step-size $\alpha$  to satisfy
\begin{align}\label{Eqn:Th_12}
\alpha \leq \frac{\left(\frac{\beta}{2L(\tau+2)+16C_1\beta}+1\right)^{\frac{1}{\tau+2}}-1}{\beta}.
\end{align}
From Lemma \ref{Lem:2}, $\Psi(z_k)$ converges linearly as (\ref{Equ:ly}).
\item Case $\tau = 0$(no delay):
Since　$~\eta_2 \leq 1$, we get $ \frac{\eta_1(1+\eta_2)^2}{\alpha} \leq 4C_1\beta$ and further obtain
\begin{align}\label{Eqn:Th_14}
\Delta_k^1 + \Delta_k^2 &= L\sum_{j=k-1}^k\|z_{j+1}-z_j\|^2+ \frac{\eta_1(1+\eta_2)^2}{\alpha}\sum_{j=k-2}^{k-1}\|z_{j+1}-z_{j}\|^2 \nonumber\\
&\leq \left(L+ 4C_1\beta\right)\sum_{j=k-2}^k\|z_{j+1}-z_j\|^2.
\end{align}
According to analogous process of case ($\tau \geq 1$), we denote $b_1=\frac{1}{4\alpha}$, $b_2=\frac{\eta_2+2\eta_2^2}{2\alpha}$, $c=L+4C_1\beta$, $k_0=2$.
From (\ref{Eqn:Th_4}), we obtain
\begin{align}\label{Eqn:Th_15}
\Psi(z_{k+1}) \leq a(\eta_1,\eta_2)\Psi(z_k) - b_1 \|z_{k+1}-z_{k} \|^2+ b_2\|z_{k}-z_{k-1} \|^2+c\sum_{j=k-k_0}^k\|z_{j+1}-z_j\|^2.
\end{align}
We know (\ref{Eqn:Th_7}) and (\ref{Eqn:Th_8}) holds for Case ($\tau=0$) as well. 
To apply Lemma \ref{Lem:2}, we need the parameters to satisfy the following inequality, combining  (\ref{Eqn:Th_7}) and (\ref{Eqn:Th_8}),
\begin{align}\label{Eqn:Th_16}
\eta_2 \frac{1+(1-C_1)\alpha\beta}{\alpha} 
&\leq b_1-\frac{c}{1-a(0,0)}\frac{1-a(0,0)^{k_0+1}}{a(0,0)^{k_0}}\nonumber\\
& = \frac{1}{4\alpha}-\frac{L+4C_1\beta}{\alpha\beta}\left((\alpha\beta+1)^{3}-1\right).
\end{align}
Thus we could choose $\eta_1,\eta_2$ from the nonempty intervals respectively,
\begin{flalign}
&\mathcal{I}_1= [0, \frac{\alpha\beta}{2}), \text{that is } 0 \leq C_1< \frac{1}{2},\label{Eqn:Th_17}\\
&\mathcal{I}_2=\left[0,\min\left\{\frac{\alpha\beta}{2}, \frac{1}{1+(1-C_1)\alpha\beta} \left(\frac{1}{4}-\frac{L+4C_1\beta}{\beta}\left((\alpha\beta+1)^{3}-1\right)
\right)
\right\}\right].\label{Eqn:Th_18}
\end{flalign}
In order to show that $\eta_2$ is selected properly, there must be
\begin{align*}
\frac{1}{4}-\frac{L+4C_1\beta}{\beta}\left((\alpha\beta+1)^{3}-1\right)>0,
\end{align*}
which is equivalent to require the step-size $\alpha$ to satisfy
\begin{align}\label{Eqn:Th_19}
\alpha \leq \frac{\left(\frac{\beta}{4L+16C_1\beta}+1\right)^{\frac{1}{3}}-1}{\beta}.
\end{align}
From lemma \ref{Lem:2},  $\Psi(z_k)$ converges linearly as (\ref{Equ:ly}).
\end{itemize}

Combining these two cases, we get the conclusions. (\ref{Equ:obj}) and (\ref{Equ:point}) follow naturally from (\ref{Equ:ly}).
\end{proof}
\end{theorem}

\begin{remark}
\begin{itemize}
	\item We can deduce from the parameter selections that convergence rate $$\rho=\frac{1+\eta_2}{1+\alpha\beta-\eta_1} < 1.$$ 
	\item If we set $\eta_1=\eta_2=0$, that is $C_1=C_2=0$, the iPIAG reduces to PIAG algorithm, and Theorem \ref{Th:1} shows similar convergent results as in paper \cite{zhang2017linear}.
\end{itemize}
\end{remark}

Moreover, we can get  two specific convergence results of PIAG-M and PIAG-NeL by setting proper inertial parameters.

\begin{corollary}\label{Corol:1}
Suppose that Assumptions \ref{Assum:1}-\ref{Assum:4} hold.
Assume the step-size and momentum parameter $\eta_1$ satisfy:
$$\alpha \leqslant \alpha_0 := 
\frac{\left[ 1 + \frac{(1-C_1)\beta}{L(\tau+1)+ C_1\beta}\right]^{\frac{1}{\tau+1}}-1}{(1-C_1)\beta},0 \leq C_1 <1,$$ 
and 
$$ \eta_1 = C_1\alpha\beta.$$
Then, the PIAG-M method converges linearly in the sense that
\begin{eqnarray*} 
\Psi(z_k) \leqslant (\frac{1}{1+\alpha\beta -\eta_1})^k \Psi(z_0), k \geq 1.
\end{eqnarray*}

Especially, when $\alpha = \alpha_0$, there is 
\begin{eqnarray}\label{Equ:es}
\Psi(z_k) \leqslant \left[1- \frac{1-C_1}{[1 +Q(\tau+1)](\tau+1)}\right]^k \Psi(z_0),
\end{eqnarray}
where $Q:= \frac{L}{\beta}$.
\begin{proof}
Front half of its proof process is analogous to Theorem \ref{Th:1} by using Lemma \ref{Lem:1} instead of Lemma \ref{Lem:2}, another slight difference is to bound $\Delta_k^1$ and $\Delta_k^2$ once again. It is trivial and omitted. 

Finally, we state the last result (\ref{Equ:es}). We replace $\alpha = \alpha_0$ and $Q = \frac{L}{\beta}$. Then we can get
\begin{eqnarray*}
(\frac{1}{1+(1-C)\alpha\beta})^k &=& \left[ 1 - \frac{1-C}{ 1 + Q(\tau+1)}\right]^{\frac{k}{1+\tau}}\\
&\leqslant&  \left[ 1-\frac{1-C}{(1 + Q(\tau+1)}\frac{1}{1+\tau} \right]^k,
\end{eqnarray*}
where the inequality is from the Bernoulli inequality $(1+ x)^r \leqslant 1+rx$ for any $x \geqslant -1$ and $0 \leqslant r \leqslant 1$.
\end{proof}
\end{corollary}


If $\eta_1=0$, this framework iPIAG reduces to PIAG-NeL.
\begin{corollary}\label{Corol:2}
Suppose that Assumptions \ref{Assum:1}-\ref{Assum:4} hold, and define a Lyapunov function
$$\Psi_2(x) := \Phi(x)-\Phi^* + \frac{1}{2\alpha}d^2(x,\mathcal{X}).$$
If the step-size and inertial parameter satisfy:
\begin{eqnarray*}
\alpha < \alpha_0 := \frac{ (W+1)^{\frac{1}{\tau+2}}-1}{\beta},
\end{eqnarray*}
where $W = \frac{\beta}{2L(\tau+2)}$, 
and
\begin{flalign*}
&\eta_2 \in \left[0,\min\left\{\frac{\alpha\beta}{2}, \frac{1}{1+\alpha\beta} \left(\frac{1}{4}-\frac{L(\tau+2)\alpha}{2\alpha\beta}\left((\alpha\beta+1)^{\tau+2}-1\right)
\right)
\right\}\right].
\end{flalign*}
Let $\bar{\mathbf{C}}:= \Psi_2(z_1)+ \frac{1+\eta_2}{1+\alpha\beta}\Psi_2(z_0)+\frac{1}{4\alpha}\|z_1-z_0\|^2_2$. Then, $\Psi_2(z_k)$ converges linearly in the sense that:
\begin{eqnarray}
\Psi_2(z_k)\leq \left(\frac{1+\eta_2}{1+\alpha\beta}\right)^{k}\bar{\mathbf{C}}, k\geq 1.
\end{eqnarray}
\begin{proof}
The proof is similar to the process to prove Theorem \ref{Th:1} and omitted. 
\end{proof}
\end{corollary}

\section{Numerical Experiments}\label{Sec:4}
This section presents two numerical examples to demonstrate the efficiency of inertial acceleration for 
PIAG. Our first example is a toy problem which was tested in \cite{aytekin2016analysis} and compares with what have been shown in their numerical part. The Second example is for Lasso problem.

\subsection{Simulation Problem}\label{Sec:4sub1}
The authors of \cite{aytekin2016analysis} solved the problem (\ref{M:Pro1}) by PIAG with setting
\begin{small}
\begin{equation*}
f_n(x) = \left \{
\begin{aligned}
&(x_n-c)^2 + \frac{1}{2}(x_{n+1} + c)^2,  n=1, \\
&\frac{1}{2}(x_{n-1} + c)^2 + \frac{1}{2}(x_{n} - c)^2,  n=N,\\
&\frac{1}{2}(x_{n-1} + c)^2 + \frac{1}{2}(x_{n} - c)^2 + \frac{1}{2}(x_{n+1} + c)^2, n=2,\cdots,N-1.
\end{aligned}
\right.
\end{equation*}
\end{small}
\begin{eqnarray*}
h(x)=\lambda_1 \|x\|_1 + I_{\mathcal{X}}(x),~~~~\mathcal{X}=\{x \geqslant 0\},
\end{eqnarray*}
for some $c \geqslant 0$. 

According to the setting of this objective, $F$ is $(N+1)$-gradient Lipschitz continuous and $F$ is $2$-stongly convex(of course satisfies quadratic growth condition with $\beta =2$). The optimizer is $x^* = \frac{\max\{0, c-\lambda_1\}}{3} \mathbf{e}_1$ with $\mathbf{e}_i$ denotes the $i$-th basis vector. As the parameter selection of paper \cite{aytekin2016analysis}, $N=100, c=3, \lambda_1=1$ and $W=4$ (each worker computes $N/W=25$ component functions' gradient). At each iteration $k$, one worker is sent to master agent with uniformly random selection and $\tau = 4$. This is only a simulation example. In practical applications, the communication between master agent and worker agents are according to the calculations' completing order rather than random choosing, and $\tau$ is unknown. Here we only use this setting to verify better performance of iPIAG than PIAG. 

\begin{figure*}[htbp]
 \centering
 \subfigure[$\alpha=0.0001$]{
  \includegraphics[width=0.47\textwidth]{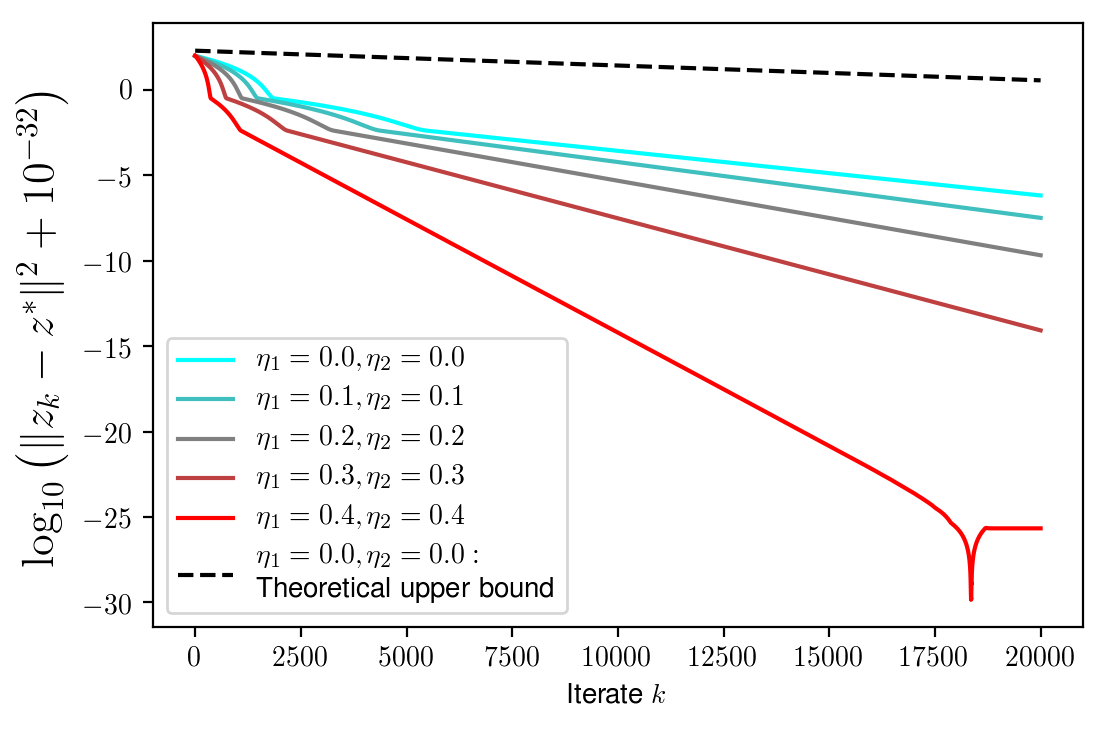}}
  \subfigure[$\alpha=0.001$]{
  \includegraphics[width=0.47\textwidth]{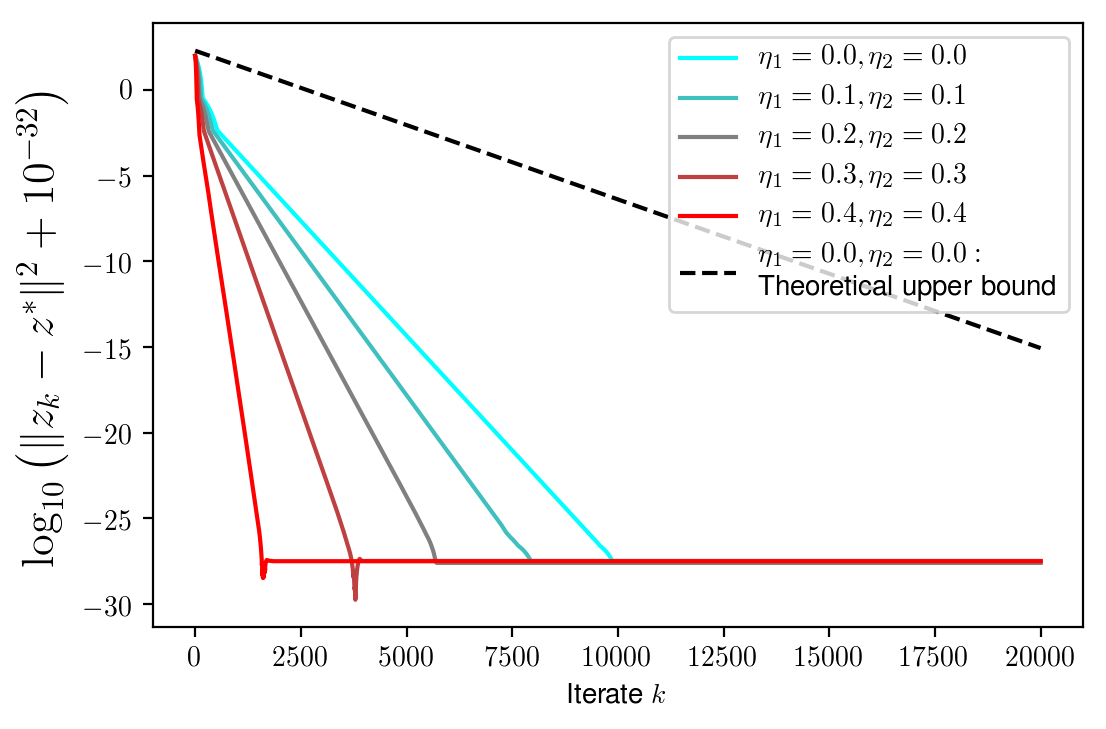}}
 \hspace{2pt}
  \subfigure[$\alpha=0.0005$]{
  \includegraphics[width=0.47\textwidth]{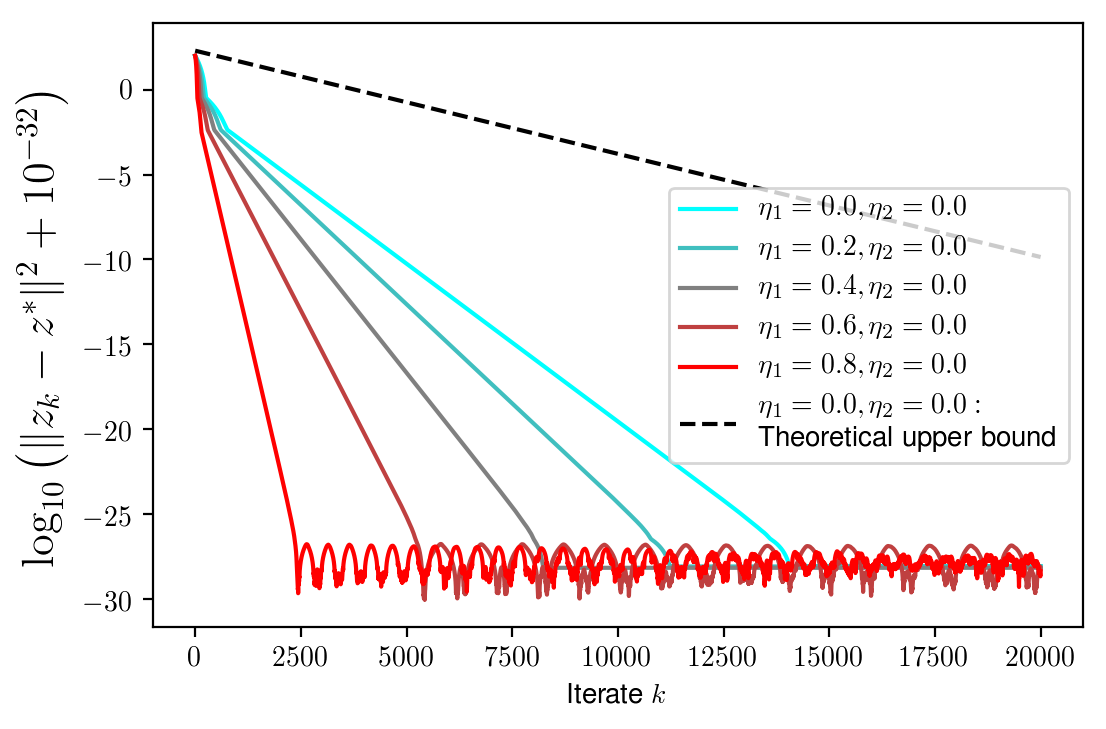}}
 \subfigure[$\alpha=0.0005$]{
  \includegraphics[width=0.47\textwidth]{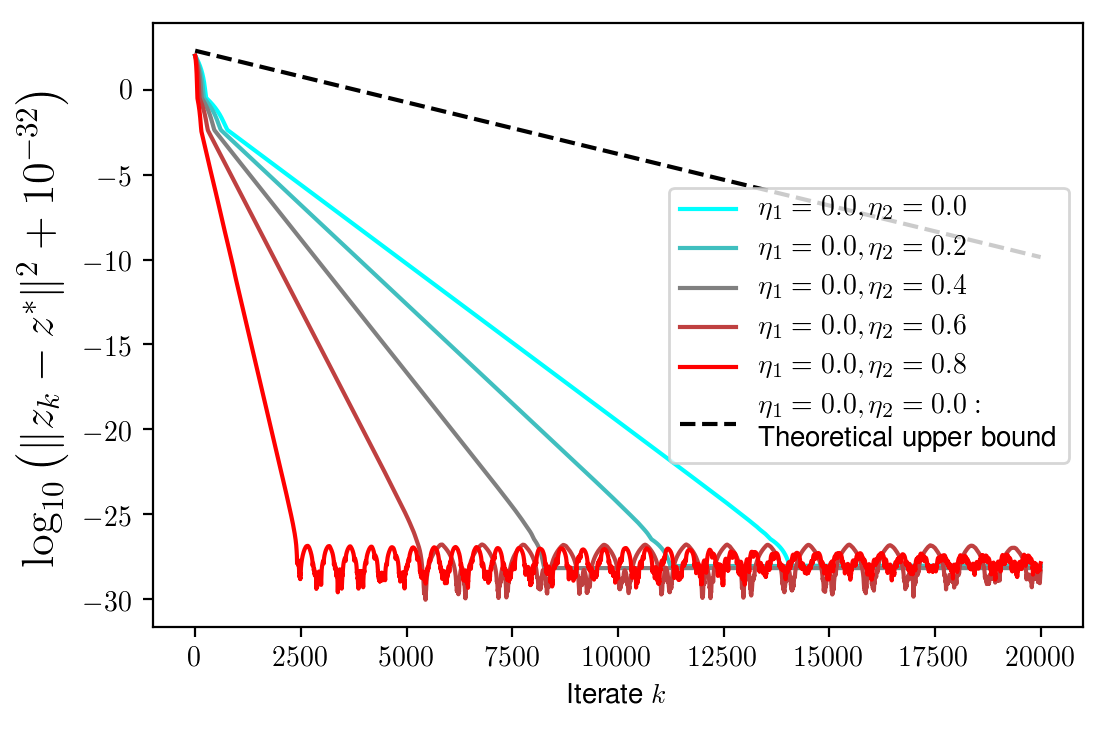}}
 \caption{Convergence of the iterates in toy problem \cite{aytekin2016analysis}. Dash-dotted line represents the best theorem upper bound in (3.16); different solid lines represent simulation results for different inertial parameter settings.\label{Fig:2}}
\end{figure*}

In figure \ref{Fig:2}, we chose the step size as the analysis result $\alpha_0$ in Theorem 1. The maximum iteration is $10000$. Figure \ref{Fig:2} show the results when running iPIAG algorithm for different step size and inertial parameter settings. We report the square of a distance from the iteration point $x_k$ to the optimal point $x^*$. 
As already mentioned, setting $\eta_1=0,\eta_2=0$, reverts the proposed iPIAG algorithm to the PAIG algorithm\cite{aytekin2016analysis}. Setting $\eta_2=0$, reduces the proposed iPIAG algorithm to the PAIG-M algorithm. Setting $\eta_2=0$, reverts the proposed iPIAG algorithm to the PAIG-NeL algorithm. 
The theoretical bound (when $\rho = \frac{1}{1 + \alpha\beta}$) is the best estimation for iPIAG and PIAG.  This test shows that increasing step size $\alpha$ could converge faster but it should not be too large since it will not converge according to our test. One can see that for a fixed step size, iPIAG algorithm outperforms PIAG algorithm.

\subsection{Lasso Problem}
In this section, we conducted an experiment on Lasso problem with asynchronous parallel system. The lasso problem to be solved is given by 
$$\text{minimize}_{~x \in \mathbb{R}^n~~} \frac{1}{2}\|Ax-b\|^2 + \lambda \|x\|_1,$$
where $A \in \mathbb{R}^{m \times n}$, $b \in \mathbb{R}^m$ and $\lambda > 0$.

In our experiment, we set that $A$ is a $300 \times 1000$ matrix whose components are random chosen from standard gaussian distribution; $x^{\ast}$ is a $1000$ variable who has $10\%$ nonzero elements; $\lambda$ is set to $\lambda = 0.2$. We tested the problem by using $3$ worker agents. All the results are average results of 100 tests.

\begin{figure*}[htbp]
  \centering
 \subfigure[$\alpha=0.00015$]{
  \includegraphics[width=0.47\textwidth]{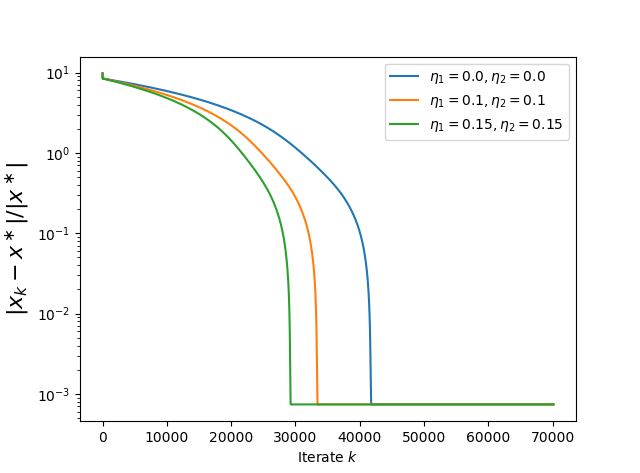}}
  \subfigure[$\alpha=0.0001$]{
  \includegraphics[width=0.47\textwidth]{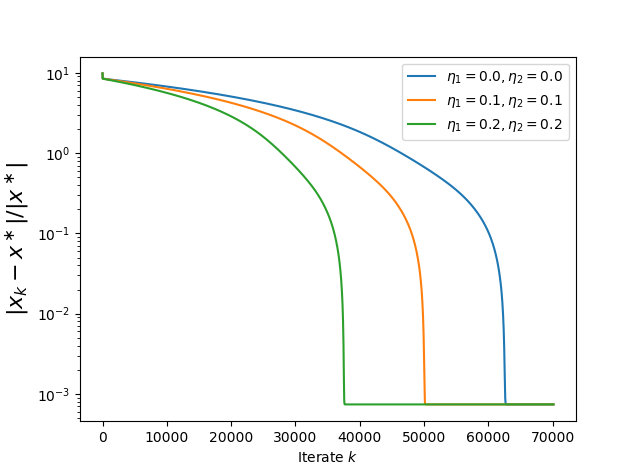}}
 \hspace{2pt}
  \subfigure[$\alpha=0.0001$]{
  \includegraphics[width=0.47\textwidth]{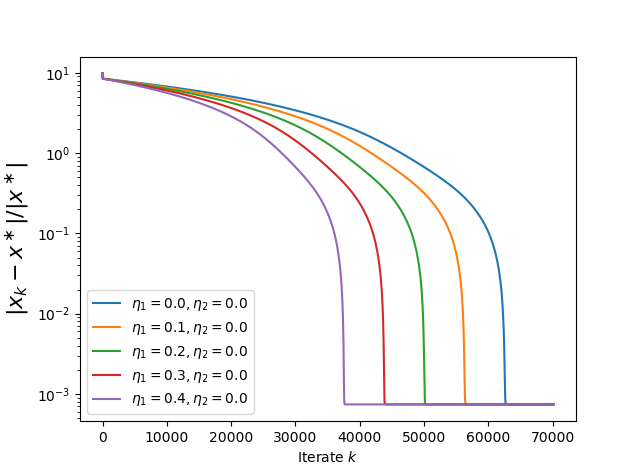}}
 \subfigure[$\alpha=0.0001$]{
  \includegraphics[width=0.47\textwidth]{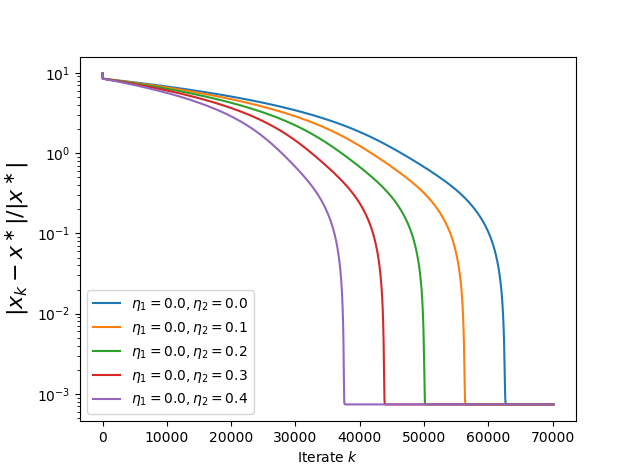}}
 \caption{Convergence of the iterates in lasso problem. Figure (a) is test whose step size is $0.00015$, the rest three cases are step size $0.0001$ for different inertial parameter settings.\label{Fig:3}}
\end{figure*}

Figure \ref{Fig:3} illustrates the convergence result with $y$-label are relative errors of iterate points to optimal point: $\|x_{k}-x^*\|/\|x^*\|$ while $x$-label are iterations $k$. From Figure \ref{Fig:3}, it can be seen that larger setting of the inertial parameters $\eta_1,\eta_2$  lead to a consistent improvement of the convergence speed.

But from the test, we observed that both of step size $\alpha$ and $\eta_1,\eta_2$ cannot be set too large, which has a consistence with the theoretical conclusions. This test also verifies that iPIAG is more powerful than PIAG.

\section{Conclusion}\label{Sec:5}
In this paper, we proposed iPIAG which is an inertial variant of the PIAG method studied in \cite{aytekin2016analysis} for solving a broad class of convex optimization problem consisting of the sum of convex gradient-Lipschitz-continuous functions and convex nonsmooth regularization function with easy to compute proximal operators.  
We proved the linear convergence of iPIAG method and provided upper bounds on the inertial and step size parameters. In particular, we showed that our results hold when the objective function  
satisfies quadratic growth condition rather than strong convexity. In our experiment, we found that proper inertial and step size parameters selection leads to significant acceleration.

\bibliographystyle{unsrt}
\bibliography{ref}

\end{document}